\documentclass{svjour3}
\journalname{MJOM}
\smartqed
\usepackage{graphicx}
\usepackage[T1]{fontenc}
\usepackage{lineno,hyperref}
\modulolinenumbers[5]

\usepackage{mathdots}       
\usepackage{algorithmic}   
\usepackage{epsfig}
\usepackage{color}
\usepackage{array}
\usepackage{anysize}
\usepackage[a4paper, total={6in, 9in}]{geometry}
\usepackage{latexsym,amssymb,amsmath,caption,setspace,titlesec}
\usepackage{mathptmx}
\usepackage{mathrsfs}
\usepackage{caption}
\usepackage{float}
\floatstyle{ruled}
\newfloat{Algorithm}{hbt}{alg}[section]
\usepackage[framed,numbered,autolinebreaks,useliterate]{mcode}
 \numberwithin{equation}{section}
\def\({\left(}
\def\){\right)}

\def\diag{\mathop{\mathrm{diag}}}

\def\sign{\mathop{\mathrm{sign}}}
\def\ind{\mathop{\mathrm{ind}}}
\def\rank{\mathop{\mathrm{rank}}}

\spnewtheorem{algorithm}{Algorithm}[section]{\bf}{\rm}

\numberwithin{equation}{section}
\numberwithin{theorem}{section}
\numberwithin{lemma}{section}
\numberwithin{corollary}{section}
\numberwithin{algorithm}{section}
\numberwithin{example}{section}
\numberwithin{proposition}{section}

\begin{document}

%
%
%
%
%
%
%
%
%

\title{Explicit Solutions of the Singular Yang--Baxter-like Matrix  Equation  and Their Numerical Computation\footnote{\bf This work has been accepted for publication in Mediterranean Journal of Mathematics and it 
will be published by April 2022.}}

\author{Ashim Kumar \and Jo\~ao R. Cardoso \and Gurjinder Singh}

\institute{Ashim Kumar \at
Department of Mathematical Sciences\\
 I.K. Gujral Punjab Technical University Jalandhar\\
 Kapurthala 144603\\  India \\
\email{ashimsingla1729@gmail.com}
\and
 Jo\~ao R. Cardoso \at 
Coimbra Polytechnic/ISEC \\
3030-199 Coimbra \\
Portugal \\
and \\
CMUC, Centre for Mathematics \\
University of Coimbra\\
3030-290 Coimbra \\
Portugal \\
\email{jocar@isec.pt}
\and
Gurjinder Singh \at
Department of Mathematical Sciences\\
 I.K. Gujral Punjab Technical University Jalandhar\\
 Kapurthala 144603\\  India \\
\email{gurjinder11@gmail.com}
}

\date{}

\maketitle


\begin{abstract}
We derive several explicit formulae for finding infinitely many solutions of the equation $AXA=XAX$, when $A$ is singular. We start by splitting the equation into a couple of linear matrix equations and then show how the projectors commuting with $A$ can be used to get families containing an infinite number of solutions. Some techniques for determining those projectors are proposed, which use, in particular, the properties of the Drazin inverse, spectral projectors, the matrix sign function, and eigenvalues. We also investigate in detail how the well-known similarity transformations like Jordan and Schur decompositions can be used to obtain new representations of the solutions. The computation of solutions by the suggested methods using finite precision arithmetic is also a concern. Difficulties arising in their implementation are identified and ideas to overcome them are discussed. Numerical experiments shed some light on the methods that may be promising for solving numerically the said matrix equation.  
\end{abstract}

\section{Introduction}

This paper deals with the equation 
\begin{equation}\label{4x}
AXA=XAX,
\end{equation}
where $A  \in \mathbb{C}^{n\times n}$ is a given complex  matrix and $X\in \mathbb{C}^{n\times n}$ has to be determined. This equation is called the {\it Yang--Baxter-like matrix equation}. If $A$ is singular (nonsingular) matrix, then the equation \eqref{4x} is said to be the  singular (nonsingular) Yang--Baxter-like matrix equation.  The equation \eqref{4x} has its origins in the classical papers by Yang \cite{Yang1} and Baxter \cite{Baxter1}. Their pioneering works have led to extensive research on the various forms of the  Yang--Baxter equation arising in braid groups, knot theory and quantum theory (see, e.g., the books \cite{Nichita,Yang2}). The YB-like equation \eqref{4x} is also known as the star-triangle-like equation in statistical mechanics; see, e.g., \cite[Part III]{McCoy}. 

A possible way of solving \eqref{4x} is to multiply out both sides, which leads to a system of $n^2$ quadratic equations with $n^2$ variables. However, this strategy may have little practical interest, unless $n$ is very small, say $n=2$ or $n=3.$ 

Note that the YB-like matrix equation \eqref{4x} has at least two trivial solutions$:$ $X = 0$ and $ X = A.$ Of course, the interest in solving it is in calculating non-trivial solutions. Discovering collections of solutions of \eqref{4x} or characterizing its full set of solutions have attracted the interest of many researchers in the last few years. Since a complete description of the solution set for an arbitrary matrix $A$ seems very challenging, many authors have been rather successful in doing so by imposing restrictive conditions on $A$. See, for instance, \cite{Cibotarica,Mansour} for $A$ idempotent, \cite{Ding15,Dong} for $A$ diagonalizable, and \cite{Tian} for matrices with rank one. 

Our interest in this paper is to solve the equation for a general singular matrix $A,$ without additional assumptions. We recall that among the published works on the YB-like equation \eqref{4x}, few are devoted to the numerical computation of its solutions. With this paper, we expect to give a contribution to fill in this gap. In our recent paper \cite{Kumar18}, we have proposed efficient and stable iterative methods for spotting commuting solutions for an arbitrary matrix $A.$ Nevertheless, those methods are not designed for determining  non-commuting solutions and there are a few cases where it is difficult to choose a good initial approximation (e.g., $A$ is non-diagonalizable). 

The principal contributions of this work w.r.t. the solutions of singular YB-like equation $AXA=XAX$ are$:$ 
\begin{enumerate}
\item[(i)] To establish a new connection between the YB-like equation and a set of two linear matrix equations, whose general solution is known; this is also valid for a nonsingular matrix  $A$  (--cf. Sect. \ref{splitting});
\item[(ii)] To explain clearly the role of  projectors commuting with $A$  in the process of deriving new families containing infinitely many solutions and to discuss how to find such projectors (--cf. Sect. \ref{proj-based});
\item [(iii)] To show how the similarity transformations can be utilized for locating more explicit representations of the solutions  (--cf. Sect. \ref{similarity});
\item[(iv)] To propose effective numerical methods for solving the singular YB-like equation, alongside with a thorough discussion of their numerical behaviour and practical clues for implementation in MATLAB (--cf. Sects. \ref{issues} and \ref{sec-experiments}).
\end{enumerate} 

By ${\mathbf 0}$ and $I,$ we mean respectively the zero and identity matrices of appropriate orders.  For a given matrix $Y,$ we denote $N(Y)$ and $R(Y)$ by the null space and the range of $Y$, respectively; $v(\lambda)$ stands for the index of a complex number $\lambda$ with respect to a square matrix $Y,$ that is,  $v(\lambda)$ is the index of the matrix $Y-\lambda I$ (check the beginning of Sect. \ref{basics} for the definition of the index of a matrix);

\section{Basics}\label{basics}

Given an arbitrary matrix $ A\in \mathbb{C}^{n\times n},$ consider the following conditions, where $X\in \mathbb{C}^{n\times n}$ is unknown$:$ 
$$\begin{array}{lll}
\textrm{(gi.1)}\ AXA=A & \quad \textrm{(gi.2)}\ XAX=X & \quad \textrm{(gi.3)}\ AX=(AX)^* \\
\textrm{(gi.4)}\ XA=(XA)^* & \quad \textrm{(gi.5)}\ AX=XA & \quad \textrm{(gi.6)}\ A^\textrm{{ind}(A)+1}X=A^{\textrm{ind}(A)},
\end{array}$$
where $Y^*$ denotes the conjugate transpose of the matrix $Y$ and $\textrm{ind}(A)$ stands for the index of a square matrix $A,$ which is the smallest non-negative integer $\ell$ such that $\rank(A^\ell)=\rank(A^{\ell+1}).$ If $m(\lambda)$ is the minimal polynomial of $A$, then $\ell$ is the multiplicity of $\lambda = 0$ as a zero of $m(\lambda)$ \cite[p. 154]{Ben}. Thus, $\ell\leq n,$ where $n$ is the order of $A.$  

A complex matrix $X \in\mathbb{C}^{n\times n}$ satisfying the condition (gi.2) is called a generalized outer inverse or a $\{2\}$-inverse of $A$, while the unique matrix $X$ verifying the conditions ({gi.1}) to ({gi.4}) is the well-known Moore--Penrose inverse of $A$, which is denoted by $A^\dag$  \cite{Penrose}; the unique matrix $X$ obeying the conditions ({gi.2}), ({gi.5}) and ({gi.6}) is the Drazin inverse, which is denoted by $A^D$ and is given by 
\begin{equation}\label{comp-drazin}
A^D=A^\ell(A^{2\ell+1})^\dag A^\ell,
\end{equation}
where $\ell \geq \ind(A).$  Other instances of generalized inverses may be defined  \cite{Ben}, but are not used in this paper. We refer the reader to \cite{Ben}, \cite[Chapter 4]{Laub}, and \cite[Section 3.6]{Lutkepohl} for the theory of generalized inverses. For both theory and computation, see \cite{Wang}. 

In the following, we revisit two important matrix decompositions, the Jordan and the Schur decompositions, whose proofs can be found in many Linear Algebra and Matrix Theory textbooks (see, for instance, \cite{Horn13}). Both decompositions will be used later in Sect.~\ref{similarity} to detect explicit solutions of the singular matrix equation $AXA=XAX.$ In addition, due to the numerical stability of the Schur decomposition, it is the basis of the algorithm that will be displayed in Figure \ref{fig1}. 

\begin{lemma} (Jordan Canonical Form) 
	Let $A\in \mathbb{C}^{n\times n}$ and let $J:=\diag(J_{n_1}(\lambda _1), \ldots,
	J_{n_s}(\lambda_s))$, ($n_1+\cdots+n_s=n$), where $ \lambda_1, \ldots,\lambda_s$ are the eigenvalues of $A$, not necessarily distinct, and 
	$J_k(\lambda)\in \mathbb{C}^{k\times k}$
	denotes a Jordan block of order $k$. Then there exists a nonsingular matrix $S\in \mathbb{C}^{n\times n}$ such that $\,A=SJS^{-1}.\,$ The Jordan matrix $J$ is unique up to the ordering of the blocks $J_k,$ but the transforming matrix $S$ is not.
\end{lemma}

For singular matrix $A$ of order $n$ with $\rank(A)=r<n,$  it is possible to reorder the Jordan blocks in a way that those blocks associated with the eigenvalue $0$ appear in the bottom-right of $J$ with decreasing size, that is, $J:=\diag\left(J_{n_1}(\lambda _1), \ldots, J_{n_p}(\lambda_p),J_{n_{p+1}}(0), \ldots, J_{n_s}(0)\right)$, with  $n_{p+1}\geq \ldots \geq n_s$ ($0\leq p\leq s$). So $A$ can be decomposed in the form
\begin{equation}\label{block1}
A=SJS^{-1}=S\, \left[\begin{array}{cc}
J_1 & {\mathbf 0} \\ {\mathbf 0} & J_0 \end{array}\right]\, S^{-1},
\end{equation}
where $J_1=\diag\left(J_{n_1}(\lambda _1), \ldots, J_{n_p}(\lambda_p)\right)$ is nonsingular and  $J_0=\diag\left(J_{n_{p+1}}(0), \ldots, J_{n_s}(0)\right)$ is nilpotent. 

\begin{lemma} (Schur Decomposition) 
	For a given matrix $A\in \mathbb{C}^{n\times n}$ there exists a unitary matrix $U$ and an upper triangular $T$ such that $\,A=UTU^{\ast},\,$ 
	where $U^{\ast}$ stands for the conjugate transpose of $U.$ The matrices $U$ and $T$ are not unique.
\end{lemma}

If $A$ is singular, then by reordering the eigenvalues in the diagonal of $T,$ where the zero eigenvalues appear in the bottom-right, the Schur decomposition of $A$ can be written in the form 
\begin{equation}\label{block3}
A=UTU^*=U\, \left[\begin{array}{cc}
B_1 & B_2 \\ {\mathbf 0} & {\mathbf 0} \end{array}\right]\,U^*,
\end{equation} 
where $B_1$ is $s \times s $ and $B_2$ is $s \times (n-s),$ with $ r=\rank(A) \leq  s \leq n-1.$ Note that $B_2$ is not, in general, the zero matrix.

Now we recall a lemma that provides an explicit solution for a well-known pair of linear matrix equations.

\begin{lemma}\label{coupled}(\cite{Cecioni,Rao})
Let $A, B, C, D\in \mathbb{C}^{n\times n}.$ The pair of matrix equations $AX=B,\ XC=D$ is consistent if and only if 
$$AD=BC,\ AA^{\dag}B=B,\ DC^\dag C=D,$$
and its general solution is given by
\begin{equation}\label{generalexpression}
X=A^\dag B+(I-A^\dag A)DC^\dag +  (I-A^\dag A)Y(I-CC^\dag),
\end{equation}
where $Y$ is an arbitrary $n\times n$ complex matrix.
\end{lemma}

Necessary and sufficient conditions for the equations $AX=B,\ XC=D$ to have a common solution are attributed to Cecioni \cite{Cecioni} and the expression \eqref{generalexpression} for a general common solution to Rao and Mitra \cite[p. 25]{Rao}. See also \cite[p. 54]{Ben} and \cite{Penrose}.

\section{Splitting the YB-Like Matrix Equation}\label{splitting}

In the next lemma, we split a general YB-like matrix equation into a system of matrix equations similar to the one in Lemma \ref{coupled}. Such a result will be useful in the next section. 
\begin{lemma}\label{connection}
	Let $A\in\mathbb{C}^{n\times n}$ be given and let $B\in\mathbb{C}^{n\times n}$ be such that the set of matrix equations
	\begin{equation}\label{set1}
	AX=B,\   XB=BA
	\end{equation}
	has at least a solution $X_0.$ Then $X_0$ is a solution of \eqref{4x}. Conversely, if $X_0$ is a solution of \eqref{4x}, then there exists a matrix $B$ such that
	$$AX_0=B,\   X_0B=BA.$$
	\end{lemma}
\begin{proof}
If $X_0$ is a solution of the simultaneous equations in \eqref{set1}, then $AX_0=B $ and $X_0B=BA.$ Therefore $\,AX_0A=BA=X_0B=X_0AX_0.\,$ Conversely, suppose that $X_0$ is a solution of $AXA=XAX$, i.e. $AX_0A=X_0AX_0.$ Letting $B:=AX_0$, we have $BA=X_0B$, which implies that $X_0$ is a solution of \eqref{set1}. \qed
\end{proof}

 Note that Lemma \ref{connection} is also valid for the nonsingular YB-like matrix equation. Using Lemmas \ref{coupled} and \ref{connection}, we must look for a matrix $B$ that makes \eqref{set1} consistent, that is, 
\begin{equation}\label{m-eq1}
ABA=B^2, \ AA^\dag B=B, \ BAB^\dag B=BA.
\end{equation}
For a given singular matrix $A$ and any of  $B$ satisfying  \eqref{m-eq1}, the matrices of the form  
\begin{equation}\label{Explicit}
X=A^\dag B+(I-A^\dag A)AB B^\dag + (I-A^\dag A)Y\big(I-BB^\dag\big),
\end{equation}
constitute an infinite family of solutions to \eqref{4x}, where $Y\in\mathbb{C}^{n\times n}$ is  arbitrary .

\section{Commuting Projectors-Based Solutions}\label{proj-based}

Discovering all the matrices $B$ in \eqref{m-eq1} may be a very hard task, apparently so difficult as solving the YB-like matrix equation. However, if $A$ is singular and $B$ is taken as in the following lemma, we have the guarantee that $B$ satisfies the conditions in \eqref{m-eq1}. Thus, many collections containing infinite solutions to the singular YB-like matrix equation can be obtained, as shown below. 

\begin{lemma}\label{define-B}
Let $A\in\mathbb{C}^{n\times n}$ be singular  and $P$ be any idempotent matrix commuting with $A$, that is, $P^2=P$ and $PA=AP.$ Then, for  $B\in\left\{A^2 P,\, A^2 (I-P)\right\}$, any matrix $X$ obtained as in  \eqref{Explicit} is a solution of  the singular YB-like matrix equation \eqref{4x}.
\end{lemma}

\begin{proof}
If $B=A^2 P$, then the equality $PA=AP$ implies that $B$ commutes with $A.$ Using the equalities $P^2=P$ and $BB^\dag B=B$, it is not difficult to show that the conditions in \eqref{m-eq1} hold for the matrix $B$ and hence the result follows. Similar arguments apply to $B=A^2(I-P).$ \qed
\end{proof}

By Lemma \ref{define-B}, we must look for matrices $P$ that are idempotent and commute with a given singular matrix $A$, in order to define $B.$ Below, several cases with examples of matrices $B$ satisfying the conditions of Lemma \ref{define-B} will be presented when $A$ is a singular matrix. 

\medskip \noindent {\bf Case 1.} $B\in\left\{{\mathbf 0},\, A^2\right\}.$

\medskip This case arises, for instance, when $P$ is a trivial commuting projector, that is, $P={\mathbf 0}$ or $P=I.$ Let us assume first that $B={\mathbf 0}.$ Now the system \eqref{set1} reduces to the matrix equation $AX={\mathbf 0}$ which is clearly solvable. From \eqref{Explicit}, its general set of solutions can be determined through the formula
 \begin{equation}\label{explicit1}
 X=(I-A^\dag A)Y. 
 \end{equation}

Geometrically speaking, the set of matrices constructed by \eqref{explicit1} is a vector subspace of $\mathbb{C}^{n\times n}$, and hence the sum of solutions of the YB-like matrix equation or a scalar multiplication yield new solutions. Since $\rank(A)=\rank(A^\dag A)$ and $\rank(I-A^\dag A)=n-\rank(A)$, such a subspace has dimension equal to $n(n-\rank(A)).$   

Now, if  we assume that $B=A^2,$ by \eqref{Explicit}, 
\begin{equation}\label{explicit3}
   X=A^\dag A^2+(I-A^\dag A)A^3 \left(A^2\right)^\dag + (I-A^\dag A)Y\big(I-A^2\left(A^2\right)^\dag\big),
   \end{equation}
   where $Y\in\mathbb{C}^{n\times n}$ is arbitrary, gives another family of solutions to \eqref{4x}. 

\medskip In the following result, we identify the solutions defined by \eqref{explicit1} and \eqref{explicit3} that commute with $A$. 

\begin{proposition}\label{pro2}
 		Let $A\in\mathbb{C}^{n\times n}$ be singular. For arbitrary matrices $Y_1,Y_2\in\mathbb{C}^{n\times n}$, the following formulae generate solutions of the equation \eqref{4x} that commute with $A$$:$
 		\begin{eqnarray}
 		X_1&=& (I-A^\dag A)Y_1(I-AA^\dag);\label{exp1}\\
 		X_2&=& A^\dag A^2+(I-A^\dag A)A^2A^\dag +(I-A^\dag A)Y_2(I-AA^\dag).\label{exp2}
 		\end{eqnarray}
 \end{proposition}
 
 \begin{proof}
	Every solution $X_1$  of $AX={\mathbf 0}$,\ $XA={\mathbf 0}$ belongs to the solution space defined by $AX={\mathbf 0}$, whose general solution is determined by  \eqref{explicit1}. Let $X_1$ be a common solution of the equations $AX={\mathbf 0}$,\ $XA={\mathbf 0}.$ Then $X_1$ commutes with $A$ and satisfies $AXA=XAX.$ By  Lemma \ref{coupled}, $X_1$ is of the form \eqref{exp1}. 
	
	Clearly, the set of matrix equations  $AX=A^2,\ XA=A^2$ is consistent and its solution set agrees with that of $AX=A^2,\ XA^2=A^3$, which is delivered by \eqref{explicit3}. Hence,  $AXA=XAX.$ If  $X_2$ is a solution of the coupled matrix equations $AX=A^2,\ XA=A^2$, then it is a commuting solution of \eqref{4x} and, again by Lemma \ref{coupled}, $X_2$ is given by \eqref{exp2}. \qed
 \end{proof}

The following lemma gives theoretical support for Case 2.

\begin{lemma}\label{Theorem21e}
Let  $A  \in \mathbb{C}^{n \times n}$ be a given singular matrix and let $M \in \mathbb{C}^{n \times n}$ be any matrix such that $AM=MA.$ Then $P_M=MM^D$ is an idempotent matrix commuting with $A.$
\end{lemma}
\begin{proof}
 Using the properties of the Drazin inverse, in particular,  $M^DMM^D=M^D$, it is easily proven that $P_M$ is an idempotent matrix.  It follows from Theorem 7 in \cite[Chapter 4]{Ben}  that $M^D$ is a polynomial in $M$ and hence $M^DA=AM^D$, because $AM=MA.$ This shows that $P_M$ commutes with $A.$  \qed
\end{proof}

 \noindent {\bf Case 2.} $B\in\left\{A^2P_M,\, A^2\(I-P_M\)\right\}.$       

\medskip It should be mentioned that the matrix $M$ in Lemma~\ref{Theorem21e} must be singular to avoid trivial cases. Examples of such matrices $M$  can be taken from the infinite collection
$$\mathfrak{M}_{\lambda_i}=\{f(A)-f(\lambda_i)I: f(x) \ \textnormal{is any polynomial  over} \ \mathbb{C}\}, $$ 
for each $ \lambda_i  \in \sigma(A)$, where $\sigma(A)$ is the spectrum of $A.$ A trivial example is to take $M=A$ yielding $P_M=AA^D$.    

\medskip The next case (Case 3) involves the matrix sign function. Before proceeding, let us recall its definition \cite[Chapter 5]{Higham08}. Let the $n\times n$ matrix $A$ have the Jordan canonical form $A=ZJZ^{-1}$ so that $J= \left[\begin{array}{cc}
	J^{(1)} & {\mathbf 0} \\ {\mathbf 0} & J^{(2)} \end{array}\right]$, where the eigenvalues of $ J^{(1)} \in \mathbb{C}^{p \times p}$ lie in  open left half-plane and those of $ J^{(2)} \in \mathbb{C}^{q \times q}$  lie in  open right half-plane. Then $S:=\sign(A)=Z\begin{bmatrix}-I_{p\times p} & 0   \\ 0 & I_{q\times q}
\end{bmatrix}Z^{-1}$ is named as the matrix sign function of $A.$ If $A$ has any eigenvalue on the imaginary axis, then $\sign(A)$ is undefined. Here $S^2=I$ and $AS=SA.$ Note also that $P=(I+ S)/2$ and $Q=(I- S)/2$ are projectors onto the invariant subspaces associated with the eigenvalues in the right half-plane and left half-plane, respectively. For more properties and approximation of the matrix sign function, see \cite{Higham08}.   

Since in this work $A$ is assumed to be singular, we cannot use directly $\sign(A)$ because it is undefined. To overcome this situation, authors in \cite{ashim} have shifted and scaled the eigenvalues of $A$ so that the matrix sign of such resulting matrices exists and commutes with $A.$ Nevertheless, there is an absence of a systematic and algorithmic approach to generate these newly matrices.  

Towards this aim, let us consider the matrix $A_\alpha :=\alpha I+A$, where $\alpha$ is a suitable complex number. Note that the scalar $\alpha$ must be carefully chosen in order to avoid the intersection of the spectrum of $A_\alpha$ with the imaginary axis.  Assuming that $A$ has at least one eigenvalue that does not lie on the imaginary axis, a simple procedure for calculating several values of $\alpha$ that leads to the acquisition of the maximal number of projectors is described as follows:
\begin{enumerate}
	\item Let $\{r_1,\ldots,r_s\}$ be the set constituted by the distinct real parts of the eigenvalues of $A$ written in ascending order, that is,  $r_1<r_2<\ldots<r_s;$
	\item For $k=1,\ldots,s-1$, choose $\alpha_k = -(r_k+r_{k+1})/2.$
\end{enumerate}

\noindent This way of calculating $\alpha_k$ guarantees that the eigenvalues of the successive $A_{\alpha_k}$ do not intersect the imaginary axis and avoids the trivial situations. That is to say, the spectrum of $A_{\alpha_k}$ does not lie entirely on either the open right half-plane or on the open left-plane, in which cases $\sign(A_{\alpha_k})=I$ or $\sign(A_{\alpha_k})=-I$. If $S_{\alpha}:=\sign(A_{\alpha})$, we see that $S_{\alpha_k}$ and $A_{\alpha_k}$ commute, because $S_{\alpha_k}$ commutes with $A_{\alpha_k}.$ Hence $\left(I+S_{\alpha_k}\right)/2$ and $\left(I-S_{\alpha_k}\right)/2$ are  projectors commuting with $A$. In the particular case when all the eigenvalues of $A$ are pure imaginary, we may consider $\widetilde{A} = -iA$ and then apply the above procedure to $\widetilde{A}$ instead of $A$.

\medskip \noindent {\bf Case 3.} $B \in \{A^2\left(\frac{I+S_{\alpha}}{2}\right), \, {A^2\left(\frac{I-S_{\alpha}}{2}\right) }\}.$	

\medskip The upcoming case depends on the spectral projectors of $A$, which have played an important role in the theory of the YB-like matrix equation, \cite{23,Ding15,spec}. Yet, there is not any definite procedure to find out them in computer algebra systems. The next proposition contributes to settle it out. 
\begin{proposition}\label{pro3}
	Let $\lambda_1,\ldots,\lambda_s$ be the distinct eigenvalues of $A \in\mathbb{C}^{n\times n}$ and assume that $G_{\lambda_i}$ denotes the spectral projector onto the generalized eigenspace $N((A-\lambda _iI)^{v(\lambda_i)})$ along $R((A-\lambda _iI)^{v(\lambda_i)})$, associated with the eigenvalue $ \lambda_i $. Then, for any $i=1,\ldots,s$, $G_{\lambda_i}$ can be represented as $G_{\lambda_i} =I-(A-\lambda_i I)(A-\lambda_i I)^D$,  where $v(\lambda_i)$ is the index of $\lambda_i.$
\end{proposition}

\begin{proof} Let $r_i:=\rank((A-\lambda_i I)^{v(\lambda_i)}).$  Since  $N((A-\lambda _iI)^{v(\lambda_i)})$ and $R((A-\lambda _iI)^{v(\lambda_i)})$ are complementary subspaces of $ \mathbb{C}^{n},$ the spectral projector $G_{\lambda_i}$ onto $N((A-\lambda _iI)^{v(\lambda_i)})$ along $R((A-\lambda _iI)^{v(\lambda_i)})$ can be written as: 	$ G_{\lambda_i}=Q_i \, \diag \left({\mathbf 0}_{r_i \times r_i},\right.$ $\, \left. I_{(n-r_i) \times (n-r_i)}\right) \, Q_i^{-1},$
	with $Q_i=[X_i | Y_i]$, in which the columns of $X_i$ and $Y_i$ are bases for $R((A-\lambda _iI)^{v(\lambda_i)})$ and $N((A-\lambda _iI)^{v(\lambda_i)})$, respectively; see, for example, \cite{Ben} and \cite[Chapters 5 and 7]{cd}.
	
	\smallskip On the other hand, the core-nilpotent decomposition of the matrix $(A-\lambda_i I)$ via $Q_i$ can be written in the form $ Q_i^{-1}(A-\lambda _iI)Q_i=\diag(C_{r_i \times r_i}, \, N_{(n-r_i) \times (n-r_i)}),
	$ 	where $C_{r_i \times r_i}$ is nonsingular, and $ N_{(n-r_i) \times (n-r_i)}$ is nilpotent of index $v(\lambda_i),$  \cite[Chapter 5, p. 397]{cd}. Now we have, $ (A-\lambda _iI)=Q_i \, \diag \left(C_{r_i \times r_i}, \, N_{(n-r_i) \times (n-r_i)}\right) \, Q_i^{-1}$ and hence the Drazin inverse of   $ (A-\lambda _iI)$ is given by $ (A-\lambda _iI)^D=Q_i \, \diag \left(C^{-1}_{r_i \times r_i}, \, {\mathbf 0}_{(n-r_i) \times (n-r_i)} \right) \, Q_i^{-1},$ \cite[Chapter 5, p. 399]{cd}. This further implies that $ 	I-(A-\lambda_i I)(A-\lambda_i I)^D=Q_i \, \diag({\mathbf 0}_{r_i \times r_i}, \, I_{(n-r_i) \times (n-r_i)}) \, Q_i^{-1},$ 	which coincides with $G_{\lambda_i}.$  \qed	
\end{proof}

Now, we revisit a well-known result, whose proof can be found in the literature (e.g., \cite{Ben,spec,cd}).  

\begin{lemma}\label{new5}
Let us assume that the notations and conditions of the Proposition \ref{pro3} are valid. Then:
\begin{itemize} 
\item[\textnormal{(a)}] $G_{\lambda_i}^2=G_{\lambda_i}$,  $AG_{\lambda_i}=G_{\lambda_i}A$, and $G_{\lambda_i}G_{\lambda_j}={\mathbf 0}$, for $i \neq j$;
\item[\textnormal{(b)}]  $P_{\lambda_i}=I-G_{\lambda_i}=(A-\lambda_i I)(A-\lambda_i I)^D$ is the complementary projector  onto $R((A-\lambda _iI)^{v(\lambda_i)})$ along $N((A-\lambda _iI)^{v(\lambda_i)})$  commuting with $A.$ In addition,  $P_{\lambda_i}P_{\lambda_j}=P_{\lambda_j}P_{\lambda_i}$;
\item[\textnormal{(c)}] $\sum_{i=1}^s G_{\lambda_i}=I$;
\item[\textnormal{(d)}] The sum of any number of matrices among the $G_{\lambda_i}$'{s} is also a commuting projector with $A.$ Thus, for any nonempty subset $\Gamma$ of $\{1,2, \dotsc, s\}$, $E_\Gamma$ is a projector commuting with $A$, where $E_\Gamma:=\sum_{i\in \Gamma} G_{\lambda_i}$;
\item[\textnormal{(e)}] $P_{\lambda_i}=\sum^s_{\substack{j=1 \\ j \neq i}} G_{\lambda_j}.$
\end{itemize}

\end{lemma}

Note that the number of projectors $E_\Gamma$'s is $2^s-1.$  Next, in Case $4$, we present a new choice of $B$  in Lemma \ref{define-B}, using the projectors described above.

\medskip\noindent {\bf Case 4.} $B \in \{A^2E_\Gamma \}.$  	      

\medskip To derive the last case (see Case 5 below), we use again the matrix sign function. It is based on the following result.  
 
\begin{proposition}\label{theoremsign}
 		Let $\lambda_1,\ldots,\lambda_s$ be the distinct eigenvalues of $A \in\mathbb{C}^{n\times n}$. For any scalar $\alpha$ and $i=1,\ldots, s$, the eigenvalues of the matrix $\hat{A}{_{\lambda_i}}=A+(\alpha-\lambda_i)G_{\lambda_i}$, consist of those of $A$, except that one eigenvalue $\lambda_i$ of $A$ is replaced by $\alpha.$ Moreover, if $\sign(\hat{A}{_{\lambda_i}})$ exists, then it commutes with $A.$ 
 \end{proposition}

\begin{proof} Let $A=P\, \diag \left(\widetilde{J}_1, \ldots, \widetilde{J}_i, \ldots, \widetilde{J}_s\right) \,  P^{-1} $ be the Jordan decomposition of $A$, where $\widetilde{J}_i$ is the Jordan segment corresponding to $\lambda_i$ and $P$ is nonsingular.  Here $ A-\lambda_i I=P \,\diag \left(\widetilde{J}_1-\lambda_i \widetilde{I}_1, \ldots, \widetilde{J}_i-\lambda_i \widetilde{I}_i, \ldots, \widetilde{J}_s-\lambda_i \widetilde{I}_s\right) \, P^{-1},$ where $\widetilde{I}_{i}$ is the identity matrix of the same order as  $\widetilde{J}_i.$ From \cite[Chapter 4, Theorem 8]{Ben}, it follows that $(A-\lambda_i I)^D=P \, \diag \left((\widetilde{J}_1-\lambda_i \widetilde{I}_1)^{-1}, \ldots,{\mathbf 0}, \ldots, (\widetilde{J}_s-\lambda_i \widetilde{I}_s)^{-1} \right) \, P^{-1}$ and we get $G_{\lambda_i} =I-(A-\lambda_i I)(A-\lambda_i I)^D=P \, \diag \left({\mathbf 0}, \ldots,\widetilde{I}_i, \right.$ $\left. \ldots, {\mathbf 0} \right) \, P^{-1}.$ Thus,  
\begin{align}\label{shift-1}
\hat{A}{_{\lambda_i}}=A+(\alpha-\lambda_i)G_{\lambda_i}&=P \, \diag \left(\widetilde{J}_1, \ldots, \widetilde{J}_i+(\alpha-\lambda_i)\widetilde{I}_i, \ldots, \widetilde{J}_s \right) \, P^{-1} \nonumber \\
&=P \, \diag \left(\widetilde{J}_1, \ldots, \widetilde{J}_i(\alpha), \ldots, \widetilde{J}_s \right) \, P^{-1},
\end{align}
where $\widetilde{J}_i(\alpha)$ is the matrix $\widetilde{J}_i$ with $\alpha$ in the place of $\lambda_i$. This shows that the eigenvalues of the matrix $\hat{A}{_{\lambda_i}}$ coincide with those of $A$ with the exception that $\lambda_i$  is replaced by  $\alpha$ in  $\hat{A}{_{\lambda_i}}.$  This proves our first claim in the proposition. \\
\indent It is clear that no eigenvalue of $\hat{A}{_{\lambda_i}}$ lies on the imaginary axis, since we are assuming that $\sign(\hat{A}{_{\lambda_i}})$ exists. Let $\hat{S}{_{\lambda_i}}=\sign(\hat{A}{_{\lambda_i}})=P \, \diag \left(\sign(\widetilde{J}_1), \ldots, \sign(\widetilde{J}_i(\alpha)), \ldots, \sign(\widetilde{J}_s) \right) \, P^{-1}$. Then a simple calculation shows that $\hat{S}{_{\lambda_i}}$ commutes with $A$ because  $\sign(\widetilde{J}_i)=\pm \widetilde{I}_{ i}$. This proves our second claim. \qed
\end{proof}

\noindent {\bf Case 5.} $B \in \{A^2\left(\frac{I+ \hat{S}{_{\lambda_i}}}{2}\right), \, A^2\left(\frac{I- \hat{S}{_{\lambda_i}}}{2}\right)\}.$

\medskip We stop here and do not pursue to attain more possibilities for $B.$ This could be considered for future works.

\section{Connections Between the Projectors and $B$}

For a given singular matrix $A$, the five cases presented in the previous section aimed at finding a commuting projector $P$ (i.e., $AP=PA$ and $P^2=P$) in order to obtain a matrix $B$ that will be inserted in \eqref{Explicit} to produce a family of solutions to the YB-like equation \eqref{4x}.

One issue arising in this approach for spotting $B$ is that distinct projectors may correspond to the same $B.$ That is to say, if $P_1$ and $P_2$ are two distinct commuting projectors then we may have $B=A^2P_1=A^2P_2$, which means that $A^2(P_1-P_2)={\mathbf 0}$, that is,  $R(P_1-P_2)\subseteq N(A^2).$ To get more insight into this connection between the projectors and $B$, we will present two simple examples.

\medskip\noindent {\bf  Example 1.} Let 
       $A=\left[\begin{array}{ccc}
       1 & 1 & 1 \\ 0 & 1 & 0 \\ 1& 1 & 1 \end{array}\right],$ which is a diagonalizable singular matrix with spectrum $\sigma(A)=\{0,1,2\}.$ 
 Solving directly the equations $AP=PA$ and $P^2=P$, we achieve a total of eight distinct commuting projectors:
 
\medskip\begin{tabular}{cccc}
 $P_1 = {\mathbf 0}$,& 
 $P_2=\frac{1}{2}\left[\begin{array}{rrr}
 	1 & 1 & 1 \\ 0 & 2 & 0 \\ 1 & -1 & 1 \end{array}\right],$ &
 $P_3=\frac{1}{2}\left[\begin{array}{rrr}
   	1 & 1 & 1 \\ 0 & 0 & 0 \\ 1 & 1 & 1 \end{array}\right],$ &
 $P_4=\left[\begin{array}{rrr}
  0 & 0 & 0 \\ 0 & 1 & 0 \\ 0 & -1 & 0 \end{array}\right],$  \\
  &&&\\
$P_5 = I$, & $P_6=I-P_2$, & $P_7=I-P_3$, & $P_8=I-P_4.$
\end{tabular}

\medskip\noindent However, there are just four distinct $B_i=A^2P_i$ ($i=1,\ldots,8$): 

\medskip\begin{tabular}{cccc}
	$B_1=A^2P_1= {\mathbf 0}$, & $B_2=A^2P_2= A^2$, & $B_3=A^2P_3$, & $B_4=A^2P_4$,\\
\end{tabular}

\medskip\noindent because $B_5=B_2$, $B_6=B_1$, $B_7=B_4$ and $B_8=B_3.$ The same four distinct $B_i$'s can be obtained by means of the sign function (Case 3) for $\alpha\in\{-5/2,-3/2,-1/2,1/2\}.$ However, Case 3 gives only six distinct projectors: $P_1,P_2,P_3,P_5,P_6,P_7,$ instead of eight projectors.  Note that the matrix sign function of $A_\alpha$ just depends on the sign of its eigenvalues, so choosing other values for $\alpha$ would not change the results. We have found those values of $\alpha$ by the method described in the previous section for Case 3. If we now find the six spectral projectors $P_{\lambda_i}$'s and $G_{\lambda_i}$'s, for all $\lambda_i\in\sigma(A)$ (see Proposition \ref{pro3} and Lemma \ref{new5}), we obtain all the commuting projectors, except the trivial ones $P_1$ and $P_5.$ Those six spectral projectors suffice to collect the four distinct matrices, $B_i$'s. 

Note that, for this matrix $A$, we can use \eqref{Explicit} to achieve four families of infinite solutions to the equation $AXA=XAX.$

\medskip\noindent {\bf  Example 2.} Let 
$A=\left[\begin{array}{ccc}
1 & 1 & 1 \\ 1 & 1 & 1 \\ 1 & 1 & 1 \end{array}\right]$, which is a diagonalizable singular matrix:
$A=S\,\diag(3,0,0)\,S^{-1},$
where $S=\left[\begin{array}{rrr}
1 & 1 & 1 \\ 1 & -1 & 1 \\ 1 & 0 & -2 \end{array}\right].$  It can be proven that all the distinct commuting projectors $P$ are given by
$$P=S\,\left[\begin{array}{cc}
\mu & {\mathbf 0} \\ {\mathbf 0} & \widetilde{P} \end{array}\right]\,S^{-1},$$ 
where $\mu\in\{0,1\}$ and $\widetilde{P}$ is any idempotent matrix of order $2.$  Since for any of those projectors $B=A^2P={\mathbf 0}$ if $\mu=0$, and $B=A^2P=A^2$ if $\mu=1$, there are just two distinct matrices:  $B={\mathbf 0}$ and $B=A^2.$ The same result is given independently by Cases 3 and 4, leading to two families of infinite solutions to the equation $AXA=XAX$ given by \eqref{Explicit}. 

\section{More Families of Explicit Solutions}\label{6}

   In this section, we provide more explicit representations for solutions to the singular YB-like equation, but now with the help of the index of $A.$
   
   \begin{proposition}\label{theorem2}
   	Assume that $A\in\mathbb{C}^{n\times n}$  is a given singular matrix such that $\ind(A)=\ell .$
   	\begin{enumerate}
   		\item[(i)] If 
   		\begin{equation}\label{Y2}
   		Y=\left(A^{\ell+1}\right)^\dag A^\ell (I-AZ)+Z,
   		\end{equation}
   		where $Z\in\mathbb{C}^{n\times n}$ is an arbitrary matrix, then, for any $V\in\mathbb{C}^{n\times n},$
   		\begin{equation}\label{formula3}
   		X=A^{\ell-1}\left(AY-I\right)V
   		\end{equation}
   		is a solution of the  YB-like matrix equation $AXA=XAX.$
   		\item[(ii)] If 
   		\begin{equation}\label{Y3}
   		Y=(I-ZA)A^\ell\left(A^{\ell+1}\right)^\dag  +Z,
   		\end{equation}
   		where $Z\in\mathbb{C}^{n\times n}$ is an arbitrary matrix, then, for any $V\in\mathbb{C}^{n\times n},$
   		\begin{equation}\label{formula4}
   		X=V\left(Y A -I\right)A^{\ell-1}
   		\end{equation}
   		is a solution of the YB-like matrix equation $AXA=XAX.$
   	\end{enumerate}
   \end{proposition}
   
   \begin{proof}
   	It is well-known that any square matrix has a Drazin inverse, which implies in particular that the matrix equation (gi.6) is solvable. From \cite[Theorem 6.3]{Laub}, it follows that  $A^{\ell+1}\left(A^{\ell+1}\right)^\dag A^\ell=A^\ell.$ Now, a simple calculation shows that the matrix $Y$ given in \eqref{Y2} is a solution of the matrix equation (gi.6), that is,  $A^{\ell+1}Y=A^\ell$, while $Y$ in \eqref{Y3} satisfies $YA^{\ell+1}=A^\ell.$ Moreover, any solution of the matrix equation $A^{\ell+1}X=A^\ell$ is of the form given in \eqref{Y2}, and any solution of  $XA^{\ell+1}=A^\ell$ can be calculated from \eqref{Y3}. The proof that both $X$ in \eqref{formula3} and $X$ in \eqref{formula4} satisfy the singular YB-like matrix equation \eqref{4x}, follows from a few matrix calculations. \qed
   \end{proof}

\section{ Solutions Based on Similarity Transformations}\label{similarity}
 

\begin{lemma}\label{lem-similar}
	Let $A,B\in \mathbb{C}^{n\times n}$ be similar matrices, that is,  $A=SBS^{-1},$ for some nonsingular complex matrix $S.$ If $Y$ is a solution of the YB-like matrix equation $BYB=YBY,$ then $X=SYS^{-1}$ is a solution of the YB-like matrix equation $AXA=XAX.$ Reciprocally, if $X$ satisfies $AXA=XAX$ then there exists $Y$ verifying $BYB=YBY$ such that $X=SYS^{-1}.$
\end{lemma}
 
The previous result, whose proof is easy, can be utilized in particular with similarity transformations like the Jordan canonical form or the Schur decomposition (--cf. Sect.~\ref{basics}).

Let us assume that $A=SJS^{-1}=S\, \left[\begin{array}{cc}
J_1 & {\mathbf 0} \\ {\mathbf 0} & J_0 \end{array}\right]\, S^{-1}$
is the Jordan decomposition of $A$, where $S$, $J_0$ and $J_1$ are as in \eqref{block1}. If 
$Y=\left[\begin{array}{cc}
	Y_{1} & Y_{2} \\ Y_{3} & Y_{4} \end{array}\right]\, $
is a solution of $YJY=JYJ,$ conformally partitioned as $J,$ then

\begin{equation}\label{sys1}
\begin{cases}
Y_{1}J_1Y_{1}+Y_{2}J_0Y_{3}=J_1Y_{1}J_1, \  \ Y_{1}J_1Y_{2}+Y_{2}J_0Y_{4}=J_1Y_{2}J_0,\\
Y_{3}J_1Y_{1}+Y_{4}J_0Y_{3}=J_0Y_{3}J_1,\ \ Y_{3}J_1Y_{2}+Y_{4}J_0Y_{4}=J_0Y_{4}J_0.
\end{cases}
\end{equation}

Hence one can determine all the solutions of equation \eqref{4x} by solving \eqref{sys1} for the matrices $Y_{\textrm{i}}$ ($i=1,2,3,4$). It turns out that building up its complete set of solutions  seems to be unattainable. However, if we consider the special case for $Y$ in which $Y_{2}={\mathbf 0}$ and $Y_{3}={\mathbf 0}$,  then \eqref{sys1} reduces to 
\begin{equation}\label{sys2}
\begin{cases}
Y_{1}J_1Y_{1}=J_1Y_{1}J_1, \\ Y_{4}J_0Y_{4}=J_0Y_{4}J_0,
\end{cases}
\end{equation}
consisting of two independent nonsingular and singular YB-like matrix equations for  $J_1$ and $J_0$, respectively. Now, we arrive at the following proposition with the help of Lemma~\ref{lem-similar}. 
\begin{proposition}\label{p1}
	Let $A\in \mathbb{C}^{n\times n}$ be a singular matrix and consider the notations used in \eqref{block1}. Then, \linebreak $ X=S\, \left[\begin{array}{cc}
	Y_1 & {\mathbf 0} \\  {\mathbf 0} & Y_4 \end{array}\right]\, S^{-1}$  is a solution of equation \eqref{4x}, where $Y_1$ and $Y_4$  satisfy their corresponding YB-like equations in \eqref{sys2}. 
\end{proposition}

Now an important issue arises: how to solve \eqref{sys2}?  A possible way is to take $Y_1=J_1$ or $Y_1={\mathbf 0}$, which satisfies the first equation in \eqref{sys2}, and then finding $Y_4$ by any of the suggested representations discussed in Sects. \ref{proj-based} and \ref{6}. Hence, a family of solutions to \eqref{4x} resulting from Proposition \ref{p1} is commuting or non-commuting according to $Y_4$ is commuting or non-commuting, respectively. 

If
$Z=\left[\begin{array}{cc}
	Z_{1} & Z_{2} \\ Z_{3} & Z_4 \end{array}\right]\, $,
which is assumed to be conformally partitioned as $T$ in \eqref{block3}, is a solution of $ZTZ=TZT$, then we come down with the next set of four equations:
\begin{equation}\label{sys4}
\begin{cases}
Z_{1}B_{1}Z_{1}+Z_{1}B_{2}Z_{3}&=B_1Z_{1}B_1+B_2Z_{3}B_1, \\ Z_{1}B_1Z_{2}+Z_{1}B_2Z_4&=B_1Z_{1}B_2+B_{2}Z_{3}B_{2},\\
Z_{3}B_1Z_{1}+Z_{3}B_2Z_{3}&= {\mathbf 0}, \\ Z_{3}B_1Z_{2}+Z_{3}B_2Z_4&= {\mathbf 0}.
\end{cases}
\end{equation}

Solving \eqref{sys4} is again a challenging task, therefore we restrict this task to the particular situation when $Z_{3}={\mathbf 0}.$  Now \eqref{sys4} becomes 
\begin{equation}\label{sys5}
\begin{cases}
Z_{1}B_{1}Z_{1}=B_1Z_{1}B_1, \\ Z_{1}B_1Z_{2}+Z_{1}B_2Z_4=B_1Z_{1}B_2,
\end{cases}
\end{equation}
which leads us to the following proposition:
\begin{proposition}\label{schur-based}
	Let $A\in \mathbb{C}^{n\times n}$ be a singular matrix of the form \eqref{block3}. Then, $ X=U\, \left[\begin{array}{cc}
	Z_1 & Z_2 \\  {\mathbf 0} & Z_4 \end{array}\right]\, U^{\ast}$, where $Z_1,$ $Z_2,$ and $Z_4$  satisfy simultaneously the equations \eqref{sys5}, is a solution of the equation \eqref{4x}. 
\end{proposition}

Some examples of solutions to (\ref{sys5}) are:
\begin{enumerate}
	\item[(i)] $Z_1={\mathbf 0}$, $Z_2$ and $Z_4$ arbitrary;
	\item[(ii)] $Z_1=B_1$, $Z_2=B_2$, and $Z_4={\mathbf 0}$;
	\item[(iii)] Any commuting solution of $Z_{1}B_{1}Z_{1}=B_1Z_{1}B_1$, along with $Z_2=B_2$ and $Z_4= {\mathbf 0}$;
	\item[(iv)] $Z_1=B_1^2B_1^D$, $Z_2=B_1B_1^DB_2$, $Z_4= {\mathbf 0}$, for the case when $B_1$ is singular.
\end{enumerate}	 

Other solutions to (\ref{sys5}) may be determined by finding $Z_1$ in the first equation $Z_{1}B_{1}Z_{1}=B_1Z_{1}B_1$, which is a YB-like equation, and then determine the unknowns $Z_2$ and $Z_4$ at a time by solving the multiple linear system 
	\begin{equation}\label{undeterm}
	\left[Z_1B_1\quad Z_1B_2\right]_{s \times n} \left[\begin{array}{c}
	Z_2 \\ Z_4
	\end{array}\right]_{n \times (n-s)}=B_1Z_1 B_2,
	\end{equation}
	provided it is consistent. For instance, if we fix $Z_1=B_1$, we know that  (\ref{undeterm}) is consistent, because $Z_2=B_2$ and $Z_4={\mathbf 0}$ satisfy it. Moreover, since $B_1$ is $s \times s $, $B_2$ is $s \times (n-s),$ with $ r=\rank(A)\leq  s \leq n-1,$ and $s<n$, it has infinitely many solutions.

\section{Numerical Issues}\label{issues}

We shall now consider the problem of solving the singular YB-like matrix equation in the finite precision environments. 

Most of the explicit formulae derived in Sect. \ref{proj-based} involve the computation of generalized inverses. We recall that the Moore--Penrose inverse is available in MATLAB through the function \texttt{pinv}, which is based on the singular value decomposition of $A.$ Many other methods and scripts are available in the literature. For instance, some iterative methods of Schulz-type (e.g., hyperpower methods) have received much attention in the last few years; see \cite{soley1} and the references therein. See also \cite{Stanimirovic}, and \cite{Wang} for the Drazin and other inverses. Formula  \eqref{Explicit} with $B$ given in Cases 3 and 5 requires the computation of the matrix sign function, which is available through many methods (check \cite[Chapter 5]{Higham08}). In Sect. \ref{sec-experiments}, a Schur decomposition-based algorithm available in \cite{mftoolbox} is used to calculate the matrix sign function. Here, the accuracy of the attained solution to singular equation \eqref{4x} depends on the difficulties arising in the intermediate estimation of those functions, viz: Moore-Penrose inverses, the sign functions, or the Drazin inverses which influence the relative error affecting the detected solutions to the singular YB-like equation.    

Although the Jordan canonical decomposition is a very important tool in the theory of matrices, we must recall that its determination using finite precision arithmetic is a very ill-conditioning problem \cite{Golub1,Kagstrom}. Excepting a few particular cases, the numerical calculation of solutions of the YB-like matrix equation by means of the Jordan decomposition must be avoided. Instead, we shall resort to the Schur decomposition, whose stability properties make it well-suited for approximations. Hence, we shall focus on designing an algorithm based on \eqref{block3}.

Even this approach is not free of risks when applied to matrices with multiple eigenvalues. We recall that the computation of repeated eigenvalues may be very sensitive to small perturbations. There are also the problems of knowing when it is reasonable to interpret a small quantity as being zero and how to correctly order the eigenvalues in the diagonal of the triangular matrix to get the form \eqref{block3}.   
 
To illustrate this, let us consider the matrix
$$A=\left[\begin{array}{rrrr}
-2 & -7 & -8 & -19\\ 
0 & -6 & -6 & -12\\
0 & 3 & 2 & 7\\
1 & 2 & 3 & 6
\end{array}\right],$$
which is nilpotent. All of its eigenvalues are zero and its Jordan canonical form is $J_4(0)$, that is, it just involves a Jordan block of order $4.$ Hence, rank$(A)=3.$ However, if we calculate the eigenvalues of $A$ in  MATLAB, which has unit roundoff $u\approx 2^{-53}$, by the function \texttt{eig}, we get 

\smallskip $\begin{array}{r}
\texttt{2.2968e-04 + 2.2974e-04i} \\
\texttt{2.2968e-04 - 2.2974e-04i}\\
\texttt{-2.2968e-04 + 2.2963e-04i}\\
\texttt{-2.2968e-04 - 2.2963e-04i}
\end{array},$

\smallskip \noindent instead of values with magnitudes more close to $u.$ This is quite expected and cannot be viewed as a failure of the algorithm used by MATLAB, because the condition number (evaluated through the function \texttt{condeig}) of the single eigenvalue of $A$ is about \texttt{4.7934e+11}. This example illustrates the shortcomings that may arise in the numerical calculation of solutions of the YB-like matrix equation by Schur decomposition when $A$ has badly conditioned eigenvalues. 

Despite  such type of examples only, the Schur decomposition performs very well for general singular matrices, as will be shown in Sect.~\ref{sec-experiments}.

\begin{figure}[ht]
	\centering
	\begin{lstlisting}
	function X = singular_yb_schur(A)
	[U,T] = schur(A,'complex');
	r = rank(T);
	n = length(A);
	E = diag(T);
	E1 = sort(abs(E));
	epsilon = E1(n-r);
	[US,TS] = ordschur(U,T,abs(E) > epsilon);
	B1 = TS(1:r,1:r); B2 =TS(1:r,r+1:n); O1 = TS(r+1:n,1:r);
	Ba = [B1^2 B1*B2];
	Bb = B1^2*B2;
	Z1 = B1;
	Xp = linsolve(Ba,Bb);
	N = null(Ba);
	D=diag(randn(1,n-r));
	Z = Xp+N(:,1:n-r)*D;
	Z2 = Z(1:r,:); Z4 = Z(r+1:n,:);
	TS1 = [Z1 Z2;O1 Z4];
	X = US*TS1*US';
	end
	\end{lstlisting}
	\caption{\small MATLAB script for finding solutions of the singular YB-like matrix equation by Schur decomposition combined with the solution of \eqref{undeterm}, with $s=r$}
	\label{fig1}
\end{figure}

In Figure \ref{fig1}, we provide a MATLAB script based on \eqref{undeterm} for obtaining solutions of the singular YB-like matrix equation. It involves the Schur decomposition, $A=UTU^\ast$, which is reordered to move all the elements in the diagonal of $T$ smaller than or equal to a certain quantity \texttt{epsilon} to the bottom-right. The tolerance \texttt{epsilon} determines what elements in the diagonal of $T$ are viewed as corresponding to the zero eigenvalue. To identify a suitable \texttt{epsilon}, we sort the eigenvalues of $T$  by increasing order of magnitude and assume that \texttt{epsilon} is the $(n-r)$-th eigenvalue in the ordered vector, where $r=\rank(A).$ Then a solution for the rank deficient linear system \eqref{undeterm} is attained by appropriate solvers.

If all of the eigenvalues of $A$ are well-conditioned or if $A$ is diagonalizable, \texttt{epsilon} is in general small; otherwise, it can be larger (say, $10^{-4}$) (--cf. Sect.~\ref{sec-experiments}). 

\section{Numerical Experiments}\label{sec-experiments}

We have considered several YB-like matrix equations corresponding to $15$ singular matrices with sizes ranging from $3\times 3$ to $20\times 20.$ The first three matrices (labelled with numbers from $1$ to $3$) are randomized and the next five matrices (from $4$ to $8$) were taken from the function \texttt{matrix} in the Matrix Computation Toolbox \cite{Higham-mct}; matrices labelled with $9$ to $15$ are academic examples, most of which are non-diagonalizable. We have selected the following four methods to get solutions of those $15$ YB-like matrix equations in MATLAB:

\begin{itemize}
	\item \texttt{alg-Case1}: script based on Case 1, with $B=A^2$, and subsequent use of \eqref{explicit3}, with $Y$ being a randomized matrix;
           \item \texttt{alg-sign}: script based on finding a $B$ as in Case 3, with $\alpha=-(r_{s-1}+r_s)/2$, and subsequent insertion in \eqref{Explicit}, with $Y$ being a randomized matrix; here  $\{r_1,\ldots,r_s\}$ is the set constituted by the distinct real parts of the eigenvalues of $A$ written in ascending order; for all the matrices in the experiments we have $s>1$;
           \item \texttt{alg-spectral}: script based on Case 4, with $B=A^2P_{\lambda_s}$, where $\lambda_s$ is the $n$-th component of the vector \texttt{eig(A)} obtained in MATLAB, and subsequent use of \eqref{Explicit}, with $Y$ being a randomized matrix;
	\item \texttt{alg-schur}: script provided in Figure \ref{fig1}.
\end{itemize}

Experiments related to other suggested formulae are not shown here. \texttt{alg-Case1}, \texttt{alg-sign}, and \linebreak \texttt{alg-spectral} involve the computation of the Moore--Penrose inverse, which has been carried out by the function \texttt{pinv} of MATLAB. The computation of the Drazin inverse in  \texttt{alg-spectral} has been based on \eqref{comp-drazin}. To estimate the quality of the approximation $\widetilde{X}$ to a solution $X$ of equation \eqref{4x}, we use the expression provided in \cite[Equation (15)]{Kumar18} for estimating the relative error, which is recalled here for convenience: 
\begin{equation}\label{rel-err}
\mathtt{est_{rel}(\widetilde{X})} = \frac{ \|R(\widetilde{X})\|}{\|M(\widetilde{X})\|\, \|\widetilde{X}\|},
\end{equation}
where $\|.\|$ stands for the Frobenius norm, $R(X):=AXA-XAX$ and $M(X):=A^T\otimes A-I\otimes(XA)-(AX)^T\otimes I \in \mathbb{C}_{n^2\times n^2}$ ($\otimes$ denotes the Kronecker product).

\begin{figure}[t]
	\centering
	\includegraphics[width=15cm]{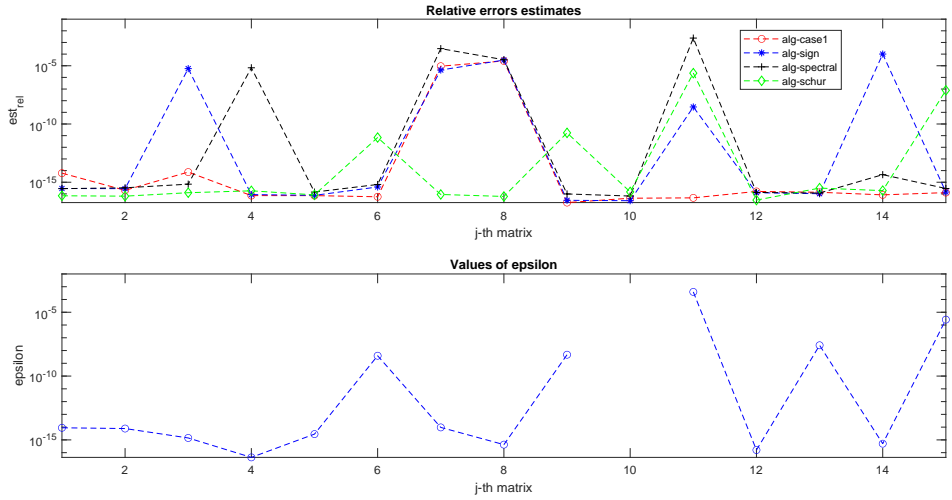}
	\caption{\small Top: relative error estimates for the solutions of \eqref{4x} obtained by \texttt{alg-Case1}, \texttt{alg-sign},  \texttt{alg-spectral},  and \texttt{alg-schur}. Bottom: values of \texttt{epsilon} defined in the script of Figure \ref{fig1}. The value of \texttt{epsilon} missing for matrix no. $10$ is exactly $zero$}
\label{fig2}
\end{figure}

At the top of Figure \ref{fig2}, we observe \texttt{alg-Case1} performs very well for all the test matrices, with the exception of matrices $7$ and $8$, where the computation of the Moore--Penrose inverses causes some difficulties. Fortunately, in these two cases, \texttt{alg-schur} gives good results. So they seem to complement very well, in the sense that when one method gives poor results the other one has a good performance. Matrices $7$ and $8$ have, respectively, sizes $19\times 19$ and $20\times 20$, and ranks $12$ and $13.$ In the case of \texttt{alg-schur}, relative errors are larger for matrices $6$, $9$, $11,$ and $15$, which are non-diagonalizable and have ill-conditioned eigenvalues as well. It is interesting to note that a comparison between both graphics shows a synchronization of the relative errors with the values of \texttt{epsilon}. 
\texttt{alg-sign} and \texttt{alg-spectral} give quite poor results for some matrices, in which large errors arise mainly in the calculation of Moore-Penrose or Drazin inverses. In the case of \texttt{alg-sign}, the choice of $\alpha$ may also influence the accuracy of the computed solutions. It is worth pointing out that arbitrary matrices $Y$ with a  large norm in \eqref{Explicit} may also cause  difficulties.

\section{Conclusions}\label{conclusions}

At this point, it is worth  highlighting the excellent features of the proposed techniques for computing solutions of singular YB-like matrix equations: 
\begin{itemize}
	\item They are valid for any singular matrix;
	\item They generate infinitely many solutions; 
	\item They perform well in finite precision environments; 
\end{itemize}
and also our main theoretical contributions:
\begin{itemize}
	\item We have provided a novel connection between the YB-like matrix equation and a well-known system of linear matrix equations, and 
	\item We have investigated the role of commuting projectors in the process of designing explicit formulae and have been able to find a large set of examples of those projectors. 
\end{itemize}

We have also overcome the main difficulties arising in the implementation of the Schur decomposition-based formula of Proposition \ref{schur-based} combined with \eqref{undeterm}, by designing an effective algorithm. We recall that many ideas of the paper (for instance, the splitting of the YB-like equation) can be extended to the nonsingular case.

\section*{Acknowledgments}
The author, Ashim Kumar, acknowledges the I. K. Gujral Punjab Technical University Jalandhar, Kapurthala for providing research support to him.

\end{document}